\documentclass[sumlimits,intlimits,reqno,10pt]{amsart}
\usepackage{amsmath, amsfonts, amssymb,amsthm}

\theoremstyle{plain}
\newtheorem{theorem}{Theorem }[section]
\newtheorem{lemma}[theorem]{Lemma }
\newtheorem{corollary}[theorem]{Corollary}
\newtheorem{proposition}[theorem]{Proposition}

\theoremstyle{definition}
\newtheorem{definition}[theorem]{Definition}
\newtheorem{remark}[theorem]{Remark}

\theoremstyle{remark}

\numberwithin{equation}{section}

\def\lr{\mathop{\longrightarrow}}\def\inlim{\mathop{\underset{t}\varprojlim}}
\def\dlim{\mathop{\underset{t}\varinjlim}}
\def\sd{semi-discrete}
\def\lc{linearly compact}

\def\hlt{Hausdorff linearly topologized}

\def\inf{\mathop{\mathrm{inf}}}\def\Ext{\mathop{\mathrm{Ext}}}
\def\Tor{\mathop{\mathrm{Tor}}}

\def\Ann{\mathop{\mathrm{Ann}}}

\def\ker{\mathop{\mathrm{ker}}}
\def\Ker{\mathop{\mathrm{ker}}}
\def\Im{\mathop{\mathrm{Im}}}
\def\Ndim{\mathop{\mathrm{Ndim}}}
\def\Soc{\mathop{\mathrm{Soc}}}\def\Ass{\mathop{\mathrm{Ass}}}
\def\Coass{\mathop{\mathrm{Coass}}}

\def\Width{\mathop{\mathrm{width}}}\def\depth{\mathop{\mathrm{depth}}}
\def\Hom{\mathop{\mathrm{Hom}}}
\def\dim{\mathop{\mathrm{dim}}}
\def\p{\mathop{\frak{p}}}
\def\m{\mathop{\frak{m}}}

\begin{document}
\title{A local homology theory  for linearly compact modules}

\author{Nguyen Tu Cuong}
\address{Nguyen Tu Cuong\\ Institute of Mathematics\\ 18 Hoang Quoc Viet Road, 10307
Hanoi, Vietnam} \email{ntcuong@math.ac.vn}
\author{Tran Tuan Nam}
\address{Tran Tuan Nam\\ The Abdus Salam International Centre for
Theoretical Physics, Trieste, Italy\\ and Ho Chi Minh University
of Pedagogy\\ 280 An Duong Vuong, District 5, Ho Chi Minh City,
Vietnam} \email{trtuannam@yahoo.com.vn}
\date{}
\maketitle

\noindent {\bf Abstract.} We introduce a local homology theory for
linearly compact modules which is in some sense dual to the local
cohomology theory of A. Grothendieck. Some basic properties  such
as the noetherianness, the vanishing and non-vanishing  of local
homology modules of linearly compact modules are proved. A duality
theory between local homology and local cohomology modules of
linearly compact modules is developed by using Matlis duality and
Macdonald duality. As consequences of the duality theorem we
obtain some generalizations of  well-known results in the theory
of local cohomology for semi-discrete linearly compact modules.
\medskip



\vskip 1.5cm

\section{Introduction}\label{S:intro}

Although the theory of local cohomology has been enveloped rapidly
for the last 40 years and proved  to be a very important tool in
algebraic geometry and commutative algebra, not so much is known
about the theory of local homology. First, E. Matlis in
\cite{matthe}, \cite{matthe2} studied the left derived functors
$L^I_\bullet(-)$ of the $I-$adic completion functor
$\Lambda_I(-)=\inlim(R/I^t\otimes_R -)$, where the ideal $I$ was
generated by a regular sequence in a local noetherian ring $R$ and
proved some duality between this functor and the local cohomology
functor by using a duality which is called today the Matlis dual
functor. Next, Simon in \cite{simsom} suggested to investigate the
module $L^I_i(M)$ when $M$ is complete with respect to the
$I$-adic topology. Later, J. P. C. Greenlees and J. P. May
\cite{greder} using the homotopy colimit, or telescope, of the
cochain of Koszul complexes to define so called local homology
groups of a module $M$ by
$$H_\bullet^I(M)=
H_\bullet(\text{Hom(Tel} K^\bullet (\underline{x}^t), M)),$$ where
$\underline{x}$ is a finitely generated system of $I$ and they
showed, under some condition on $\underline x$ which are
automatically satisfied when $R$ is noetherian, that the left
derived functors $L^I_\bullet(-)$ of the $I-$adic completion can
be computed in terms of these local homology groups. Then came
the work of L. Alonso Tarr\' io, A. Jeremias L\' opez and J.
Lipman \cite{aloloc}, they gave in that paper  a  sheafified
derived-category generalization of Greenlees-May results for a
quasi-compact separated scheme.  Note that a strong connection
between local cohomology and local homology was shown  in  \cite{aloloc} and  \cite{greder}. Recently in
\cite{cuothe}, we defined  the $i$-th {\it local homology module}
$H^I_i(M)$ of an $R-$module $M$ with respect to the ideal $I$ by
$$H^I_i(M)=\underset{t}{\underleftarrow{\lim}}
{\Tor}^R_i(R/I^t,M).$$  We also proved
in \cite{cuothe}  many basic properties of local homology modules
and that $H^I_i(M)\cong L^I_i(M)$ when $M$ is artinian. Hence we
can say that there exists  a theory
for the left derived functors $L^I_\bullet(-)$ of the $I-$adic
completion (as the local homology functors) on the category of
artinian modules over  noetherian local rings parallel to the theory of local cohomology
functors on the category of noetherian modules. However, while the
local cohomology functors $H_I^\bullet(-)$ are still defined as
the right derived functors of the $I$-torsion functor $\Gamma
_I(-)= \underset{t}{\underrightarrow{\lim}}\Hom _R(R/I^t,- )$ for
not finitely generated modules, our definition of local homology
module above may not coincide  with $L^I_i(M)$ in this case. One of the most important 
reasons  is that, even if the ring $R$ is noetherian, the $I$-adic completion functor
$\Lambda_I(-)$ is neither left nor right exact on the category of
all $R$-modules. Fortunately, it was shown by results of C. U.
Jensen in \cite{jenles} that the inverse limit functors and
therefore the local homology functors still have good behaviour on
the category of linearly compact modules. The purpose of this
paper is towards a local homology  theory for  linearly compact
modules. It should be mentioned that the concept of linearly
compact spaces was first introduced by Lefschetz \cite{lefalg} for
vector spaces of infinite dimension and it was then generalized
for modules by D. Zelinsky \cite{zellin} and  I. G. Macdonald
\cite{macdua}. It was also studied by other authors such as H.
Leptin \cite{leflin}, C. U. Jensen \cite{jenles}, H. Z\"oschinger
\cite{zoslin} $\ldots.$ The class of linearly compact modules is
very large, it contains many important classes of modules such as
the class of artinian modules, or the class of finitely generated
modules over a complete ring.

The organization of our paper is as follows. In section
\ref{S:mdcptt}  we recall the concepts of linearly compact and semi-discrete linearly compact modules by using the terminology of Macdonald
\cite{macdua} and their basic facts . For any $R-$module $N$ and a \lc\  $R-$module $M$
we show that there exists uniquely a topology induced by a free
resolution of $N$ for $\Ext_R^i(N,M)$, and in addition $N$ is
finitely generated, for $\Tor^R_i(N,M)$; moreover these modules are
\lc.

In section \ref{S:mddddp} we present some basic properties of
local homology modules of \lc\ modules such as the local homology
functor $H^I_i(-)$ is closed in the category of \lc\ modules
(Proposition \ref{P:dddpcpttlcptt}).
Proposition \ref{P:mdttdnvgrl} shows that our definition of local
homology modules can be identified with the definition of local
homology modules of J. P. C.
Greenlees and J. P. May \cite[2.4]{greder} in the category of
linearly compact modules. 

In section \ref{S:ttvkttdddp} we study the vanishing and
non-vanishing of local homology modules.   Let $M$ be a   linearly
compact $R-$module with $\Ndim M=d,$ then $H^I_i (M)=0$ for all
$i> d$ (Theorem \ref{T:dlttdddpcp}). It was proved in \cite[4.8,
4.10]{cuothe} that $\Ndim M= \max \{ i\mid H^{\m}_i(M)\not=0\}$ if
$M$ is an artinian module over a local ring $(R,\m)$, where $\Ndim
M$ is the noetherian dimension defined by N. R. Roberts
\cite{robkru} (see also \cite{kirdim}). Unfortunately, as in a
personal communication of H. Z\"oschinger, he gave us the
existence of semi-discrete \lc\ modules $K$ of noetherian
dimension 1 such that  $H^{\m}_i(K)=0$ for all non-negative
integers $i$.  However, we can prove in Theorem
\ref{T:dlttdddpcpnrr} that the above equality still holds for
semi-discrete \lc\ modules with $\Ndim M\not = 1,$ moreover $\Ndim
\Gamma_{\m}(M) = \max \big\{ i\mid H^{\m}_i(M)\not=0\big\}$ if
$\Gamma_{\m}(M)\not= 0.$ 

In section \ref{S:tcntedddp} we show that local homology modules
$H^{\m}_i(M)$ of a semi-discrete linearly compact module $M$ over
a noetherian local ring $(R,\m)$ are noetherian modules on the
$\m-$adic completion $\widehat{R}$ of $R$ (Theorem  \ref{T:dltcnotemddddp}) On the other hand, for
any ideal $I$,  $H^I_d(M)$ is a noetherian $\Lambda_I(R)-$module
provided $M$ is a semi-discrete linearly compact $R-$module with
the Noetherian dimension $\Ndim M=d$ (Theorem
\ref{T:dlmddddpnotetcd}).

Section \ref{S:dnmdn} is devoted to study duality. In this section
$(R,\m)$ is a noetherian local ring and the topology on $R$ is the
$\m-$adic topology. Let  $E(R/\m)$ be the injective envelope of
$R/\m$ and $M$  a \hlt\ $R-$module. Then the {\it Macdonald dual}
$M^*$ of $M$ is defined by $M^*=Hom(M,E(R/\m))$ the set of
continuous homomorphisms of $R-$modules. Note by Macdonald
\cite[5.8]{macdua} that $M$ is a semi-discrete module if and only
if $D(M)=M^*$, where $D(M)=\Hom(M,E(R/\m))$ is the Matlis dual of
$M$. The main result of this section is Theorem
\ref{T:dldndddpvddddp} which gives a duality between local
cohomology modules and local homology modules. 

In the last section, based on  the duality theorem
\ref{T:dldndddpvddddp}  and the properties of local homology
modules in previous sections  we  can extend  some well-known
properties of local cohomology of finitely generated modules for
semi-discrete \lc\ modules.  

In this paper, the terminology "isomorphism" means "algebraic
isomorphism" and  "topological isomorphism" means
"algebraic isomorphism with the homomorphisms (and its inverse)
are continuous".

\medskip
\section{ Linearly compact modules}\label{S:mdcptt}
\medskip

First we recall the concept of  linearly compact modules by using
the terminology of I. G. Macdonald  \cite{macdua}  and some their
basic properties. Let $M$ be a topological $R-$module. A {\it
nucleus} of $M$ is a neighbourhood of the zero element of $M,$ and
a {\it nuclear base} of $M$ is a base for the nuclei of $M.$ If
$N$ is a submodule of $M$ which contains a nucleus then $N$ is
open (and therefore closed) in $M$ and $M/N$ is discrete. $M$ is
Hausdorff if and only if the intersection of all the nuclei of $M$
is $0.$ $M$ is said to be {\it linearly topologized} if $M$ has a
nuclear base $\mathcal{M}$ consisting of submodules.

\begin{definition}\label{D:dnmdcptt} A Hausdorff linearly topologized $R-$module
$M$ is said to be {\it linearly compact} if $M$ has the following
property: if $\mathcal{F}$ is a family of closed cosets (i.e.,
cosets of closed submodules) in $M$ which has the finite
intersection property, then the cosets in $\mathcal{F}$ have a
non-empty intersection.
\end{definition}

It should be noted that an artinian $R-$module is linearly compact
with the discrete topology (see \cite[3.10]{macdua}).

\begin{remark} \label{R:rmcsmdtt} Let $M$ be an $R-$module. If $\mathcal{M}$ is a
family of submodules of $M$ satisfying the conditions:

\noindent (i) For all $N_1,\ N_2 \in \mathcal{M}$ there is an
$N_3\in \mathcal{M}$ such that $N_3\subseteq N_1\cap N_2,$

\noindent (ii) For an element $x\in M$ and $N\in \mathcal{M}$
there is a nucleus $U$ of $R$ such that $Ux \subseteq N,$

\noindent then $\mathcal{M}$ is a base of a linear topology on $M$
(see \cite[2.1]{macdua}).
\end{remark}

The following  properties of linearly compact modules are often
used in this paper.

\begin{lemma} \label{L:tcmdcptt4} {\rm  (see \cite[\S 3]{macdua})}  (i) Let $M$ be
a \hlt\ $R-$module, $N$ a closed submodule of $M.$ Then $M$ is
linearly compact if and only if $N$ and $M/N$ are linearly
compact.

\noindent (ii)  Let $f: M\longrightarrow N$ be a continuous
homomorphism of \hlt\ $R-$modules. If $M$ is linearly compact,
then $f(M)$ is linearly compact and therefore $f$ is a closed map.

\noindent (iii)   If $\{M_i\}_{i\in I}$ is a family of linearly
compact $R-$modules, then $\underset{i\in I}\prod M_i$ is linearly
compact with the product topology.

\noindent (iv)   The inverse limit of a system of linearly compact
$R-$modules and continuous homomorphisms is linearly compact with
the obvious topology.
\end{lemma}

\begin{lemma} \label{L:lnkhmdcptt} {\rm (see \cite[7.1]{jenles})} Let $\{M_t\}$ be
an inverse system of linearly compact modules with continuous
homomorphisms. Then $\underset{t}{\underleftarrow{\lim}^i} M_t =0$
for all $i>0.$ Therefore, if  $$0\longrightarrow
\{M_t\}\longrightarrow\{N_t\} \longrightarrow \{P_t\}
\longrightarrow 0$$ is
 a short exact sequence of inverse systems of
$R-$modules, then  the sequence of inverse limits
$$0\longrightarrow \underset{t}{\underleftarrow{\lim}} M_t
\longrightarrow \underset{t}{\underleftarrow{\lim}} N_t
\longrightarrow \underset{t}{\underleftarrow{\lim}} P_t
\longrightarrow 0$$is exact.
\end{lemma}

Let $M$ be a linearly compact $R-$module and $F$ a  free
$R-$module with a base $\big\{ e_i\big\}_{i\in J}.$ We can define
the topology on $\Hom_R(F,M)$ as the product topology via the
isomorphism $\Hom_R(F,M)\cong M^J,$ where $M^J=\underset{i\in
J}\prod M_i$ with $M_i=M$  for all $i\in J.$ Then $\Hom_R(F,M)$ is
a \lc\ $R-$module by \ref{L:tcmdcptt4} (iii). Moreover, if $h: F
\longrightarrow F'$ is a homomorphism of free $R-$modules, the
induced homomorphism $h^*: \Hom_R(F',M)\longrightarrow
\Hom_R(F,M)$ is continuous by \cite[7.4]{jenles}. Let now
$${\bold F_{\bullet}}: \ \ldots\longrightarrow F_i \longrightarrow \ldots \longrightarrow F_1\longrightarrow F_0
\longrightarrow N \longrightarrow 0.$$ a free resolution  of an
$R$-module $N$. Then $\Ext_R^i(N,M)$ is a linearly topologized
$R$-module with the quotient topology  of $\Hom(F_i, M)$. This
topology on $\Ext_R^i(N,M)$ is called the topology induced by the
free resolution ${\bold F_{\bullet}}$ of $N$.

\begin{lemma} \label{L:extcptt}   Let  $M$ be a linearly compact $R-$module
and $N$  an   $R-$module. Then for all $i\geq 0,\ \Ext_R^i(N,M)$
is a linearly compact $R-$module with the topology induced  by a
free resolution of $N$ and this topology is independent  of the choice  of free
resolutions of $N.$  Moreover, if $f: N\longrightarrow N'$ is a
homomorphism of   $R-$modules, then the induced homomorphism $
\Ext_R^i(N',M) \longrightarrow \Ext_R^i(N,M)$ is continuous.
\end{lemma}
\begin{proof} Let ${\bold F_{\bullet}}$  be a free resolution of $N.$
It follows as about  that
 $\Hom_R({\bold F_{\bullet}},M)$ is a complex of linearly compact modules
 with continuous homomorphisms.
Therefore $\Ext^i_R(N,M)=H^i(\Hom_R({\bold F_{\bullet}},M))$ is
\lc \   by \ref{L:tcmdcptt4} (i), (ii).
Let now   ${\bold G_{\bullet}}$ be  a second free resolution of
$N$.
Then we get  a quasi-isomorphism of complexes $\varphi_{\bullet} :
{\bold F_{\bullet}} \longrightarrow {\bold G_{\bullet}}$ lifting
the identity map of $N$.  Therefore the induced homomorphism
$$\bar{\varphi}_i: H^i({\Hom}_R({\bold F_{\bullet}},M))\longrightarrow H^i({\Hom}_R({\bold G_{\bullet}},M))$$
is a topological isomorphism by \cite[7.4]{jenles} and
\ref{L:tcmdcptt4} (i), (ii) for all $i$. Similarly we can prove
for the rest statement \end{proof}

Let $N$ be a finitely generated $R-$module and  $${\bold
F_{\bullet}}=\ \ldots\longrightarrow F_i \longrightarrow \ldots
\longrightarrow F_1\longrightarrow F_0 \longrightarrow N
\longrightarrow 0$$ a free resolution of $N$ with the finitely
generated free modules. As above, we can define for a linearly
compact module $M$ a topology on $\Tor^R_i(N,M)$ induced from the
product topology of $F_i\otimes_R M$. Then by an argument
analogous to that used for the proof of Lemma \ref{L:extcptt}, we
get the following lemma.

\begin{lemma}\label{L:torcptt}    Let $N$ be a finitely generated $R-$module
and $M$ a linearly compact $R-$module. Then  $\Tor^R_i(N,M)$ is a
linearly compact $R-$module with the topology induced  by a free
resolution of $N$ (consisting of finitely generated free modules)
and this topology is independent of the choice of free resolutions of $N.$
Moreover, if $f: N\longrightarrow N'$ is a homomorphism of
finitely generated $R-$modules, then the induced homomorphism
$\psi_{i,M} : \Tor^R_i(N,M) \longrightarrow \Tor^R_i(N',M)$ is
continuous.\end{lemma}

The next result is often used in the sequel.

\begin{lemma}\label{L:ghnghtor}  Let $N$ be a finitely generated $R-$module
and $\{M_t\}$ an inverse system of linearly compact $R-$modules
with continuous homomorphisms. Then for all $i\geq 0,$ $\{\Tor^R_i
(N, M_t )\}$ forms an inverse system of linearly compact modules
with continuous homomorphisms.  Moreover, we have
$${\Tor}^R_i (N, \underset{t}{\underleftarrow{\lim}} M_t )
\cong \underset{t}{\underleftarrow{\lim}} {\Tor}^R_i (N, M_t ).$$
\end{lemma}
\begin{proof} Let ${\bold F_\bullet}$ be a free resolution of $N$
with finitely generated free $R-$modules. Since $\{M_t\}$ is an
inverse system of linearly compact modules with continuous
homomorphisms, $\{F_i\otimes_R M_t\}$ forms an inverse system of
linearly compact modules with continuous homomorphisms for all
$i\geq 0$ by \ref{L:tcmdcptt4} (iii). Then $ \{\Tor^R_i (N, M_t
)\}$ forms an inverse system of linearly compact modules with
continuous homomorphisms. Moreover
$$ {\bold F_\bullet} \otimes_R \underset{t}{\underleftarrow{\lim}} M_t
\cong  \underset{t}{\underleftarrow{\lim}}({\bold F_\bullet}
\otimes_R  M_t ),$$since the inverse limit commutes with the
direct product and
$$H_i (\underset{t}{\underleftarrow{\lim}}({\bold F_\bullet} \otimes_R M_t ))
\cong\underset{t}{\underleftarrow{\lim}} H_i ({\bold F_\bullet}
\otimes_R M_t )$$by \ref{L:lnkhmdcptt} and \cite[6.1, Theorem
1]{norani}. This finishes the proof. \end{proof}

A \hlt\ $R-$module $M$ is called {\it semi-discrete} if every
submodule of $M$ is closed. Thus a discrete $R-$module is
semi-discrete. The class of semi-discrete linearly compact modules
contains all artinian modules. Moreover, it also contains all
finitely generated modules in  case $R$ is a complete local
noetherian ring (see \cite[7.3]{macdua}). It should be mentioned
here that our notions of linearly compact and semi-discrete
modules follow Macdonald's definitions in \cite{macdua}. Therefore
the notion of linearly compact modules defined by H. Z\"oschinger
in \cite{zoslin} is different to our notion of linearly compact
modules, but it is coincident with the terminology of
semi-discrete linearly compact modules in this paper.

Denote by $L(M)$ the sum of all artinian submodules of $M,$ we
have the following  properties of semi-discrete linearly compact
modules.

\begin{lemma}\label{L:lmcpttnrrat}  {\rm (see \cite[1 (L5)]{zoslin})}  Let $M$ be
a semi-discrete linearly compact $R-$module. Then $L(M)$ is an
artinian module.
\end{lemma}

We now recall the concept of {\it co-associated primes} of a
module (see \cite{chacop}, \cite{yascoa}, \cite{zoslin}). A prime
ideal $\p$ is called {\it co-associated} to a non-zero $R-$module
$M$ if there is an artinian homomorphic image $L$ of $M$ with
$\p=\Ann_R L.$ The set of all co-associated primes to $M$ is
denoted by $\Coass_R (M).$ $M$ is called $\p-${\it coprimary} if
$\Coass_R(M)=\big\{ \p\big\}.$ A module is called {\it
sum-irreducible} if it can not be written as a sum of two proper
submodules.  A sum-irreducible module $M$  is $\p-$coprimary,
where $\p=\big\{ x\in R/ x M\not= M \big\}$  (see
\cite[2]{chacop}).

\begin{lemma}\label{L:cpttnrrcahh} {\rm (see \cite[1 (L3,L4)]{zoslin})}  Let $M$
be a semi-discrete linearly compact $R-$module. Then $M$ can be
written as a finite sum of sum-irreducible modules and therefore
the set $\Coass(M)$ is finite.
\end{lemma}
\medskip

\section{ Local homology modules of linearly compact
modules}\label{S:mddddp}
\medskip

Let $I$ be an ideal of $R,$  the {\it i-th local homology} module
$H^I_i (M)$ of an $R-$module $M$ with respect to $I$ is defined by
(see \cite[3.1]{cuothe})
$$H^I_i (M) = \underset{t}{\underleftarrow{\lim}}{\Tor}^R_i (R/I^t , M).$$It is clear that $H^I_0(M)\cong \Lambda_I(M),$
in which $\Lambda_I(M)=\underset{t}{\underleftarrow{\lim}} M/I^tM$
the $I-$adic completion of $M.$

\begin{remark}\label{R:rmdndddp} (i) As $I^t{\Tor}^R_i(R/I^t,M)=0,$
${\Tor}^R_i(M/I^tM,N)$ has a natural  structure as a module over
the ring $R/I^t$ for all  $t>0.$ Then $H^I_i(M) =
\underset{t}{\underleftarrow{\lim}}\Tor^R_i (R/I^t , M)$ has a
natural structure as a module over the ring $\Lambda_I(R)=\inlim
R/I^t.$

\noindent (ii) If $M$ is a finitely generated $R-$module, then
$H^I_i(M)=0$ for all $i>0$ (see \cite[3.2 (ii)]{cuothe}).
\end{remark}

\begin{lemma}\label{L:dddpvpksvdn}  {\rm (see \cite[\S 3]{cuothe})} Let $I$ be an
ideal generated by elements $ x_1, x_2,\ldots, x_r$  and
$H_i(\underline{x}(t),M)$ the $i-$th Koszul homology module of $M$
with respect to the sequence $\underline{x}(t)=(x^t_1,\ldots,
x^t_r).$ Then for all $i\geq 0,$

\noindent (i) $H^I_i (M) \cong \underset{t}{\underleftarrow{\lim}}
H_i (\underline{x}(t), M),$

\noindent(ii)   $H^I_i(M)$ is $I-$separated, it means that
$\underset{t>0}\bigcap I^t H^I_i(M)=0.$\end{lemma}

Let $M$ be a \lc\ $R-$module. Then $\Tor^R_i(R/I^t,M)$ is also a
\lc\ $R-$module by the topology defined as in \ref{L:torcptt}, so
we have an induced topology on the local homology module
$H^I_i(M).$

\begin{proposition}\label{P:dddpcpttlcptt}  Let $M$ be a linearly compact
$R-$module. Then for all $i\geq 0,$ $H^I_i(M)$ is  a linearly
compact $R-$module.
\end{proposition}
\begin{proof} It follows from \ref{L:torcptt} that $\{\Tor^R_i(R/I^t,M)\}_t$
forms an inverse system of linearly compact modules with
continuous homomorphisms. Hence  $H^I_i(M)$ is also a linearly
compact $R-$module by \ref{L:tcmdcptt4} (iv). \end{proof}

The following proposition shows that local homology modules can be
commuted  with inverse limits of inverse systems of linearly
compact $R-$modules with continuous homomorphisms.

\begin{proposition}\label{P:dddpghghn}  Let  $\{ M_s\}$ be  an inverse system of
linearly compact $R-$modules with the continuous homomorphisms.
Then
$$H^I_i(\underset{s}\varprojlim
M_s)\cong\underset{s}\varprojlim H^I_i( M_s).$$
\end{proposition}
\begin{proof} Note that inverse limits are commuted. Therefore
\begin{align*} H^I_i(\underset{s}\varprojlim M_s)&=\underset{t}\varprojlim {\Tor}^R_i(R/I^t,\underset{s}\varprojlim M_s)\\
&\cong \underset{t}\varprojlim\underset{s}\varprojlim {\Tor}^R_i(R/I^t, M_s)\\
&\cong \underset{s}\varprojlim\underset{t}\varprojlim
{\Tor}^R_i(R/I^t, M_s) =\underset{s}\varprojlim H^I_i(
M_s)\end{align*} by \ref{L:ghnghtor}.\end{proof}

Let $L^I_i$ be the $i-$th left derived functor of the $I-$adic
completion functor $\Lambda_I.$ The next result
shows that in case $M$ is linearly compact, the local homology
module $H^I_i (M)$ is isomorphic to the module $L^I_i(M),$ thus
our definition of local homology modules can be identified with the
definition of J. P. C. Greenlees and J. P. May   (see
\cite[2.4]{greder}).

\begin{proposition}\label{P:mdttdnvgrl}  Let  $M$ be  a linearly compact
$R-$module. Then $$ H^I_i (M)\cong  L_i^I (M)$$ for all $i\geq
0.$\end{proposition}
\begin{proof} For all $i\geq 0$ we have a
short exact sequence by \cite[1.1]{greder}, $$0\longrightarrow
\underset{t}{\underleftarrow{\lim}^1}{\Tor}^R_{i+1}(R/I^t,M)\longrightarrow
L^I_i(M)\longrightarrow  H^I_i(M)\longrightarrow 0.$$Moreover, it
follows from \ref{L:torcptt} that $\{\Tor^R_{i+1}(R/I^t,M)\}$
forms   an inverse system of  linearly compact modules with
continuous homomorphisms.  Hence, by \ref{L:lnkhmdcptt}
$$ \underset{t}{\underleftarrow{\lim}^1}{\Tor}^R_{i+1} (R/I^t,M)=0$$  and the conclusion follows.
\end{proof}

The following corollary is an immediate consequence of
\ref{P:mdttdnvgrl}.

\begin{corollary}\label{C:hqdkdddp}  Let
$$0 \longrightarrow M' \longrightarrow M \longrightarrow M" \longrightarrow 0$$ be a short exact sequence of linearly compact modules. Then we have a long exact sequence of local homology modules
$$\cdots \longrightarrow H^I_i (M') \longrightarrow H^I_i (M)\longrightarrow H^I_i (M")\longrightarrow \ \ $$
$$\cdots \longrightarrow H^I_0 (M') \longrightarrow H^I_0 (M)\longrightarrow H^I_0 (M")\longrightarrow 0.$$\end{corollary}

The following theorem gives us a characterization of $I-$separated
modules.

\begin{theorem}\label{T:dldtmditach} Let $M$ be a linearly compact $R-$module.
The following statements  are equivalent:

\noindent(i) $M$ is $I-$separated, it means that
$\underset{t>0}\bigcap I^t M =0.$

\noindent(ii) $M$ is complete with respect to the $I-$adic
topology, it means that $\Lambda_I(M) \cong M.$

\noindent(iii)  $H^I_0 (M) \cong M,\ H^I_i (M)= 0$ for all $i>0.$
\end{theorem}

To prove Theorem \ref{T:dldtmditach}, we need the two auxiliary
lemmas. The first lemma shows that local homology modules $H^I_i
(M)$ are $\Lambda_I-$acyclic for all $i>0.$

\begin{lemma}\label{L:dddpacylic}  Let  $M$ be a linearly compact $R-$module.
Then for all $j\geq 0,$
$$H^I_i (H^I_j (M)) \cong \begin{cases} H^I_j (M),& i= 0,  \\
0, & i>0.\end{cases}$$\end{lemma}
\begin{proof} It follows from \ref{L:torcptt} that $\big\{\Tor^R_j(R/I^t,
M)\big\}_t$ forms an inverse system of \lc\ $R-$modules with the
continuous homomorphisms. Then we have by \ref{P:dddpghghn} and
\ref{L:dddpvpksvdn} (i), \begin{align*}  H^I_i (H^I_j (M))
&= H^I_i(\underset{t}{\underleftarrow{\lim}}{\Tor}^R_j (R/I^t , M))\\
&\cong \underset{t}{\underleftarrow{\lim}} H^I_i({\Tor}^R_j (R/I^t , M)) \\
&\cong
\underset{t}{\underleftarrow{\lim}}\underset{s}{\underleftarrow{\lim}}
H_i (\underline{x}(s), {\Tor}^R_j (R/I^t , M)),
\end{align*} in which $\underline{x} = (x_1 , \ldots ,x_r )$ is a system of generators of $I$ and $\underline{x}(s) = (x_1^s , \ldots ,x_r^s).$
Since $\underline{x}(s){\Tor}^R_j (R/I^t , M) = 0$ for all $s\geq
t,$ we get
$$\underset{s}{\underleftarrow{\lim}} H_i (\underline{x}(s), {\Tor}^R_j (R/I^t , M)) \cong
\begin{cases} {\Tor}^R_j (R/I^t , M),& i=0,\\
0,& i>0.\end{cases}$$This finishes the proof.
\end{proof}

\begin{lemma}\label{L:dcddgiaomd}  Let  $M$ be a linearly compact $R-$module.
Then
$$H^I_i (\underset{t>0}\bigcap I^t M) \cong \begin{cases} 0,& i=0,\\
H^I_i (M),& i>0.\end{cases}$$\end{lemma} \begin{proof} From the
short exact sequence of linearly compact $R-$modules
$$0\longrightarrow I^tM \longrightarrow M\ \longrightarrow M/I^tM
\longrightarrow 0$$  for all $t>0$ we derive  by
\ref{L:lnkhmdcptt} a short exact sequence of  linearly compact
$R-$modules
$$0\longrightarrow \underset{t>0}\bigcap I^tM \longrightarrow M \longrightarrow \Lambda_I(M) \longrightarrow 0.$$
Hence we get a long exact sequence of local homology modules
$$\cdots \longrightarrow H^I_{i+1} (\Lambda_I(M)) \longrightarrow H^I_i (\underset{t>0}\bigcap I^tM) \longrightarrow H^I_i (M) \longrightarrow   H^I_i (\Lambda_I(M))\longrightarrow \ $$
$$\cdots \longrightarrow H^I_1 (\Lambda_I(M)) \longrightarrow H^I_0 (\underset{t>0}\bigcap I^tM) \longrightarrow H^I_0 (M) \longrightarrow   H^I_0 (\Lambda_I(M))\longrightarrow 0. $$
The lemma now follows from \ref{L:dddpacylic}. \end{proof}

\begin{proof}[Proof of Theorem \ref{T:dldtmditach}] $(ii)\Leftrightarrow
(i)$ is clear from the short exact sequence
$$0\longrightarrow \underset{t>0}\bigcap I^t M \longrightarrow M\longrightarrow \Lambda_I(M)\longrightarrow 0.$$
\noindent $(i) \Rightarrow (iii).$  We have  $H^I_0
(M)\cong\Lambda_I(M) \cong M.$ Combining  \ref{L:dcddgiaomd} with
(i) gives  $H^I_i (M)\cong H^I_i (\underset{t>0}\bigcap I^t M) =0$
for all $i > 0.$

\noindent $(iii) \Rightarrow (ii)$ is trivial.
\end{proof}

From Theorem \ref{T:dldtmditach} we have the following criterion
for a finitely generated module over a local noetherian ring to be
linearly compact.

\begin{corollary}\label{C:dkmdhhscpttvdd}   Let $(R,\m)$ be a local noetherian
ring and $M$ a finitely generated $R-$module.  Then $M$ is a \lc\
$R-$module if and only if $M$ is complete with respect to the $\m-$adic
topology.
\end{corollary}
\begin{proof}  Since $M$ is a finitely generated $R-$module, $M$ is
$\m-$separated. Thus,  if $M$ is a \lc\ $R-$module,
$\Lambda_{\m}(M)\cong M$ by \ref{T:dldtmditach}. Conversely, if
$M$ is complete in $\m-$adic topology, we have $M\cong \inlim
M/\m^t M.$ Therefore $M$ is a \lc\ $R-$module by \ref{L:tcmdcptt4}
(iv), as $M/\m^t M$ are artinian $R-$modules for all $t>0.$
\end{proof}
\medskip

\section{ Vanishing and non-vanishing of local homology
modules}\label{S:ttvkttdddp}
\medskip

Recall  that $L(M)$ is the sum of all artianian submodules of $M$
and $\Soc(M)$ the socle of $M$ is  the sum of all simple
submodules of $M$. The $I-$torsion functor $\Gamma_I$  is defined
by $\Gamma _I(M) = \underset{t>0}\cup ({0:_M} {I^t}).$ To prove
the vanishing and non-vanishing theorems of local cohomology
modules, we need the following lemmas.

\begin{lemma}\label{L:dktdptdcqdddp}  Let $M$ be a semi-discrete linearly compact
$R-$module. Then $H^I_0(M)=0$ if and only if $xM=M$ for some $x\in
I.$\end{lemma} \begin{proof} By \cite[2.5]{cuothe}, $H^I_0(M)=0$
if and only if $IM=M.$  Hence the result follows from
\ref{L:cpttnrrcahh}  and \cite[2.9]{chacop}.\end{proof}

\begin{lemma}\label{L:dddpmdsb0tt}  Let $M$ be a semi-discrete linearly compact
$R-$module and $\Soc(M)=0$. Then $$H^I_i(M)=0$$for all
$i>0.$
\end{lemma}
\begin{proof} Combining \ref{L:dcddgiaomd} with \ref{L:dktdptdcqdddp}we may assume, by replacing $M$ with
$\underset{t>0}\bigcap I^tM,$  that there is an $x\in
I$ such that $x M=M.$ As $\Soc(M)=0,$ it follows from \cite[1.6
(b)]{zoslin} that $0:_Mx=0.$ Thus we have an isomorphism
$M\overset{x}\cong M.$ It induces an isomorphism
$$ H^I_i(M)\overset{x}\cong H^I_i(M)$$for all $i>0.$ By \ref{L:dddpvpksvdn} (ii), we have
$$ H^I_i(M) = x  H^I_i(M) = \underset{t>0}\bigcap x^t  H^I_i(M)=0$$for all $i>0.$
\end{proof}

\begin{lemma}\label{L:lmbtcmdx}  Let $M$ be a semi-discrete linearly compact
$R-$module. Then there are only finitely many distinct  maximal
ideals  ${\m}_1, {\m}_2, \ldots, {\m}_n$ of $R$ such that
$$L(M) = \underset{j=1}{\overset{n}\bigoplus} \Gamma_{\m_j}(M).$$
\end{lemma}

\begin{proof} By \ref{L:lmcpttnrrat}, $L(M)$ is an artinian $R-$module. Thus, by
virtue of \cite[1.4]{shaame}  there are  finitely many distinct
maximal ideals ${\m}_1, {\m}_2, \ldots, {\m}_n$ of $R$ such that
$$L(M) = \underset{j=1}{\overset{n}\bigoplus}
\Gamma_{\m_j}(L(M))\subseteq \underset{j=1}{\overset{n}\bigoplus}
\Gamma_{\m_j}(M).$$ Therefore it  remains to show that $
\Gamma_{\m}(M)$ is artinian for any maximal ideal $\m$ of $R$.
Indeed,   there is from \cite[Theorem]{zoslin}  a short exact
sequence $0 \longrightarrow N \longrightarrow M \longrightarrow A
\longrightarrow 0$, where $N$ is finitely generated and $A$ is
artinian.  Then  we have an exact sequence
$$0 \longrightarrow \Gamma_{\m}(N) \longrightarrow \Gamma_{\m}(M)
\longrightarrow \Gamma_{\m}(A).$$ Obviously, $\Gamma_{\m}(A)$ is
an artinian $R-$module, $\Gamma_{\m}(N)$ is a finitely generated
$R-$module annihilated by a power of $\m$, and hence it is of
finite length. So $\Gamma_{\m}(M)$ is an artinian $R-$module as
required.
\end{proof}

\begin{lemma}\label{L:dddpbtdddpmdmx}  Let $M$ be a semi-discrete linearly compact
$R-$module. Then there are only finitely many distinct  maximal
ideals ${\m}_1, {\m}_2, \ldots, {\m}_n$ of $R$ such that
$$H^I_i(M)\cong \underset{j=1}{\overset{n}\bigoplus} H^I_i(\Gamma_{\m_j}(M))$$for all $ i>0,$ and the following sequence is exact
$$0\longrightarrow \underset{j=1}{\overset{n}\bigoplus} H^I_0(\Gamma_{\m_j}(M))\longrightarrow H^I_0(M)\longrightarrow
H^I_0(M/\underset{j=1}{\overset{n}\bigoplus}\Gamma_{\m_j}(M))\longrightarrow
0.$$\end{lemma}

\begin{proof}  The short exact sequence of linearly
compact $R-$modules$$0\longrightarrow L(M)\longrightarrow
M\longrightarrow M/L(M)\longrightarrow 0$$gives rise  to a long
exact sequence of local homology modules
$$\ldots\longrightarrow H^I_{i+1}(M/L(M))\longrightarrow H^I_i(L(M))\longrightarrow H^I_i(M)\longrightarrow H^I_i(M/L(M))\longrightarrow\ldots.$$
By \ref{L:dddpmdsb0tt},  $H^I_i(M/L(M))=  0$ for all $i>0,$ as
$\Soc(M/L(M))=0.$ Then we get $H^I_i(M)\cong H^I_i(L(M))$ for all
$i>0$ and the short exact sequence
$$0\longrightarrow H^I_0(L(M))\longrightarrow H^I_0(M)\longrightarrow H^I_0(M/L(M))\longrightarrow 0.$$
Now the conclusion follows from \ref{L:lmbtcmdx}.\end{proof}

We have an immediate consequence of \ref{L:lmbtcmdx} and
\ref{L:dddpbtdddpmdmx} for the local case.

\begin{corollary}\label{C:dddpmbdddpgamamvdp}  Let $(R,\m)$ be a local noetherian ring
and $M$ a semi-discrete linearly compact $R-$module. Then
$$L(M)= \Gamma_{\m}(M), \ \ H^I_i(M)\cong  H^I_i(\Gamma_{\m}(M))$$for all $ i>0,$ and the following sequence is exact
$$0\longrightarrow H^I_0(\Gamma_{\m}(M))\longrightarrow H^I_0(M)\longrightarrow
H^I_0(M/\Gamma_{\m}(M))\longrightarrow 0.$$
\end{corollary}

We now recall the  concept of {\it Noetherian dimension} of an
$R-$module $M$ denoted by $\Ndim M.$ Note that the notion of
Noetherian dimension was introduced first by R. N. Roberts
\cite{robkru} by the name Krull dimension. Later, D. Kirby
\cite{kirdim} changed this terminology of Roberts and refereed to
{\it Noetherian dimension} to avoid confusion with well-know Krull
dimension of finitely generated modules.  Let $M$ be an
$R-$module. When $M=0$ we put $\Ndim M = -1.$ Then by induction,
for any ordinal $\alpha,$ we put $\Ndim M = \alpha$ when (i)
$\Ndim M < \alpha$ is false, and (ii) for every ascending chain
$M_0 \subseteq M_1 \subseteq \ldots$ of submodules of $M,$ there
exists a positive integer $m_0 $  such that $\Ndim(M_{m+1} /M_m )<
\alpha$ for all $m \geq m_0$. Thus $M$ is non-zero and finitely
generated if and only if $\Ndim M = 0.$ If $0 \longrightarrow M"
\longrightarrow M\longrightarrow M'\longrightarrow 0$ is a short
exact sequence of $R-$modules, then $\Ndim M= \max\{ \Ndim M",
\Ndim M'\}.$

\begin{remark}\label{R:cntrmarcpnrr}   (i) In case $M$ is an artinian module, $\Ndim M
<\infty$ (see \cite{robkru}). More general, if $M$ is a
semi-discrete \lc\ module, there is a short exact sequence $0
\longrightarrow N \longrightarrow M \longrightarrow A
\longrightarrow 0$ where $N$ is finitely generated and $A$ is
artinian (see \cite[Theorem]{zoslin})). Hence $\Ndim M = \max\{
\Ndim N, \Ndim A\} < \infty.$

(ii) If $M$ is an artinian $R-$module or more general, a
semi-discrete \lc\ $R-$module,  then $\Ndim M \leqslant \max\{\dim
R/\p \mid \p \in \Coass(M)\}.$ Especially,  if  $M$ is an artinian
module over a   complete local noetherian ring $(R,\m),$ $\Ndim M
= \max\{\dim R/\p \mid \p \in \Coass(M)\}$ (see
\cite[2.10]{yasmag}).
\end{remark}

\begin{lemma}\label{L:cntmoochiamx} Let $M$ be an $R-$module with $\Ndim M=d>0$
and $x\in R$ such that $x M=M.$ Then $$\Ndim 0:_M x\leqslant
d-1.$$
\end{lemma}
\begin{proof} Consider the ascending chain
$$0\subseteq 0:_Mx\subseteq 0:_Mx^2 \subseteq \ldots.$$As $\Ndim M=d,$ there exists a positive integer $n$ such that
$\Ndim (0:_Mx^{n+1}/0:_Mx^n) \leqslant  d-1.$ Since $x M=M,$  the
homomorphism $0:_Mx^{n+1}/0:_Mx^n \overset{x^n}\longrightarrow
0:_Mx$ is an isomorphism. Therefore $\Ndim 0:_M x\leqslant d-1.$
\end{proof}

\begin{theorem}\label{T:dlttdddpcp}   Let $M$ be a   linearly compact $R-$module
with $\Ndim M=d.$ Then
$$H^I_i (M)=0$$
for all $i> d.$\end{theorem} \begin{proof}  Let $\mathcal{M}$ be a
nuclear base of $M$. Then, by \cite[3.11]{macdua},  $M =
\underset{U\in\mathcal{M}}{\underleftarrow{\lim}} M/U.$ It follows
from \ref{P:dddpghghn} that
$$H^I_i(M) \cong  \underset{U\in\mathcal M}{\underleftarrow{\lim}}
H^I_i(M/U).$$ Note that  $M/U$ is a discrete linearly compact
$R-$module with $\Ndim M/U \leqslant \Ndim M.$ Thus we only need
to prove the theorem  for the case $M$ is a discrete linearly
compact $R-$module. Let $L(M)$ be the sum of all artinian
$R-$submodules of $M,$ by \ref{L:lmcpttnrrat}, $L(M)$ is atinian.
From the proof of \ref{L:dddpbtdddpmdmx}, we have the isomorphisms
$$H^I_i(M)\cong H^I_i(L(M))$$ for all $ i>0$. As $\Ndim L(M)\leqslant
\Ndim M=d,$ $H^I_i(L(M)) = 0$ for all $i>d$ by  \cite[4.8]{cuothe}
and then the proof is complete. \end{proof}

\begin{remark}\label{R:rmdlttcpvdzg}   In \cite[4.8, 4.10]{cuothe} we proved that if
$M$ is an  artinian module on a local noetherian ring $(R, \m),$
then
$$\Ndim M= \max \big\{ i\mid H^{\m}_i(M)\not=0\big\},$$
where we use the convention that $\max(\emptyset)=-1$. Therefore
it raises to the following  natural question that whether the
above equality  holds true when $M$ is a semi-discrete linearly
compact module? Unfortunately, the answer is negative in general.
The following counter-example is due to H. Z\" oschinger. Let
$(R,\m)$ be a complete local  noetherian domain of dimension $1$
and $K$ the field of fractions of $R.$  Consider $K$ as an
$R-$module. Then $\Soc(K)=0$ and $\Coass(K)=\{0\},$  therefore
$\Ndim K =1$ by \cite[1.6 (a)]{zoslin}.  Since $K/R$ is artinian,
it follows by \cite[Theorem]{zoslin} that  $K$ is a semi-discrete
linearly compact $R-$module. As $xK=K$ for any non-zero element $
x\in \m$, $H^{\m}_0(K)=0$ by \ref{L:dktdptdcqdddp}. Moreover, we
obtain by \ref{L:dddpmdsb0tt} that $H^{\m}_i(K)=0$ for all $i>0$ .
Thus
$$\Ndim K =1\not = -1=  \max \big\{ i\mid H^{\m}_i(K)\not=0\big\}.$$
\end{remark}

However, the following theorem gives an affirmative answer for the question when $\Ndim M \not= 1$.

\begin{theorem}\label{T:dlttdddpcpnrr} Let $(R,\m)$ be a local noetherian  ring
and $M$ a non zero semi-discrete linearly compact $R-$module. Then

\noindent (i)\ \ $\Ndim \Gamma_{\m}(M) = \max \big\{ i\mid
H^{\m}_i(M)\not=0\big\}$  if $ \Gamma_{\m}(M)\not= 0;$

\noindent (ii)  $\Ndim M= \max \big\{ i\mid
H^{\m}_i(M)\not=0\big\}$ if   $\Ndim M\not= 1.$
\end{theorem}
\begin{proof} (i)  Since $\Gamma_{\m}(M)$ is the artinian $R-$module,
we obtain from \cite[4.8, 4.10]{cuothe} that  $$\Ndim
\Gamma_{\m}(M) = \max \big\{ i\mid
H^{\m}_i(\Gamma_{\m}(M))\not=0\big\}.$$ Thus (i) follows from
\ref{C:dddpmbdddpgamamvdp}.

\noindent (ii) First, note by virtue of \cite[1.6 (a)]{zoslin} and
\ref{R:cntrmarcpnrr} (ii) that if $\Soc (M) = 0$ then $\Ndim M
\leqslant 1$. If $\Gamma_{\m}(M)=0$ then $\Soc (M)=0$ by
\ref{L:lmbtcmdx}. So we get from the hypothesis that $\Ndim M =0$.
It follows that $M$  is a finitely generated $R-$module and
$H^{\m}_0(M)\cong \widehat M \not = 0$, where $\widehat M$ is the
$\m-$adic completion of $M$. Thus (ii) is proved in this case.
Assume now that $\Gamma_{\m}(M)\not = 0$. By (i) we have only to
show that $\Ndim M=\Ndim \Gamma_{\m}(M).$ Indeed, it is trivial
for the case $\Ndim M=0$. Let $\Ndim M>1.$ From the short exact
sequence $0 \longrightarrow \Gamma_{\m}(M) \longrightarrow M
\longrightarrow M/\Gamma_{\m}(M) \longrightarrow 0$ we get
$$\Ndim M = \max\{\Ndim \Gamma_{\m}(M), \Ndim M/\Gamma_{\m}(M)\}.$$  Since $\Soc(M/\Gamma_{\m}(M))=0$,  $\Ndim
(M/\Gamma_{\m}(M))\leqslant 1$.  Thus $\Ndim M=\Ndim
\Gamma_{\m}(M)$ as required. \end{proof}

A sequence of elements $x_1 ,\ldots ,x_r$ in $R$ is said to be an
{\it $M-$coregular} sequence (see \cite[3.1]{ooimat}) if $0:_M
(x_1 , \ldots , x_r ) \not= 0$ and $0:_M (x_1 , \ldots , x_{i-1} )
\overset{x_i}\longrightarrow  0:_M (x_1 , \ldots , x_{i-1} ) $ is
surjective for $i=1,\ldots , r.$  We denote by $\Width_I(M)$ the
supremum of the lengths of all  maximal $M-$coregular sequences in
the ideal $I.$  Note by \ref{R:cntrmarcpnrr} (i) and
\ref{L:cntmoochiamx}   that
$${\Width}_I(M)\leqslant \Ndim M<\infty$$  when $M$ is a
semi-discrete \lc\ $R-$module.

\begin{theorem}\label{T:dtdrdddpcpnrr}  Let $M$ be a semi-discrete linearly
compact $R-$module and $I$ an ideal of $R$ such that
$0:_MI\not=0.$ Then all maximal $M-$coregular sequences in $I$
have the same length. Moreover
$${\Width}_I(M)=\inf\{i/H^I_i(M)\not=0\}.$$\end{theorem}
\begin{proof} It is sufficient to prove that if  $\{x_1,
x_2,\ldots,x_n\}$  is  a maximal  $M-$coregular sequence in $ I$,
then $H^I_i (M)=0$ for all $i<n,$ and  $H^I_n (M)\not=0.$  We
argue by the induction on $n$.  When $n=0,$ there does not exists
an element $x$ in $I$ such that $xM=M.$ Then $H^I_0(M)\not=0$ by
\ref{L:dktdptdcqdddp}.

Let $n>0.$ The short exact sequence $$0\longrightarrow 0:_Mx_1
\longrightarrow M\overset{x_1}\longrightarrow M\longrightarrow 0$$
gives rise  to a long exact sequence
$$\ldots\longrightarrow H^I_i(M)\overset{x_1}\longrightarrow H^I_i(M) \longrightarrow H^I_{i-1}(0:_M{x_1})\longrightarrow\ldots.$$
By the inductive hypothesis,  $H^I_i(0:_M{x_1})=0 $ for all
$i<n-1$ and $H^I_{n-1}(0:_M{x_1})\not=0.$  Therefore by virtue of
\ref{L:dddpvpksvdn} (ii),
$H^I_i(M)=x_1H^I_i(M)=\underset{t>0}\bigcap x_1^tH^I_i(M)=0$ for
all $i<n.$  Now, it follows from  the exact sequence
$$\ldots\longrightarrow H^I_n(M)\overset{x_1}\longrightarrow H^I_n(M) \longrightarrow H^I_{n-1}(0:_M{x_1})\longrightarrow 0$$
and  $H^I_{n-1}(0:_M{x_1})\not=0$ that $H^I_n(M)\not=0$ as required.
\end{proof}

\begin{remark}\label{R:rmdkgamanotoct}   We give here an example which shows that the
condition $\Gamma _{\m} (M)\not= 0$ in Theorem
\ref{T:dlttdddpcpnrr} (i) is needful. Let $R$ be the ring and  $K$
the $R-$module as in \ref{R:rmdlttcpvdzg}. Set $M=N\oplus K$,
where $N$ is a finitely generated $R-$module satisfying $\depth
_{\m} M \geq 1$. Then, it is easy to check that $\Gamma _{\m} (M)=
0$ and $H_i^{\m} (M)\cong H_i^{\m} (K) = 0$ for all $i\geq 1$ and
$H_0^{\m} (M)\cong H_0^{\m} (N)\cong \widehat N\not= 0$. Therefore
$$\Ndim \Gamma_{\m}(M) = -1 \not= 0=\max \big\{ i\mid
H^{\m}_i(M)\not=0\big\}.$$
\end{remark}

We have seen in Remark \ref{R:rmdlttcpvdzg} the existence of a
non-zero semi-discrete \lc\ module $K$ such that $H^{\m}_i (K)=0$
for all $i\geq 0$. Below, we give a characterization for this
class of semi-discrete \lc\ modules. This corollary also shows
that we can not drop the condition $0:_M I\not= 0$ in the
assumption of Theorem \ref{T:dtdrdddpcpnrr}.

\begin{corollary}\label{C:hqddpttmi}  Let $(R,\m)$ be a local noetherian ring
and $M$  a non-zero semi-discrete \lc\ module. Then $H^{\m}_i(M)
=0$ for all $i\geq 0$ if
and only if  there exists an element $x\in \m$ such that $xM=M$ and $0:_Mx=0$.
\end{corollary}
\begin{proof} Let $H^{\m}_i(M) =0$ for all $i\geq 0$. We obtain by
\ref{L:dktdptdcqdddp} that $xM=M$ for some $x\in \m$. On the other
hand, it follows from  the short exact sequence $0\longrightarrow
0:_Mx \longrightarrow M \overset{x}\longrightarrow M
\longrightarrow 0$ that $H^{\m}_i(0:_Mx)=0$ for all $i\geq 0.$
Since $0:_Mx$ is artinian by \cite[Corollary 1 (b0)]{zoslin},
$0:_Mx=0$ by \cite[4.10]{cuothe}. Conversely, suppose that $xM=M$
and $0:_Mx=0$, then for all $i\geq 0$ $$H^{\m}_i(M) = x
H^{\m}_i(M) = \underset{t>0}\bigcap x^tH^{\m}_i(M) =0$$
by \ref{L:dddpvpksvdn} (ii).
\end{proof}
\medskip

\section{ Noetherian local homology modules}\label{S:tcntedddp}
\medskip

First, the following criterion for a module to
be noetherian is useful for the investigation of the noetherian property of local
homology modules.

\begin{lemma}\label{L:tcnote}   Let $J$ be a finitely generated  ideal of a
commutative  ring $R$ such that $R$ is complete with respect to the $J-$adic
topology and $M$ an $R-$module. If $M/JM$ is a  noetherian
$R-$module and $M$ is $J-$separated (i. e., $\underset{t>0}\bigcap
J^tM =0$), then $M$ is a noetherian $R-$module.\end{lemma}
\begin{proof} Set
$$K = \underset{t\geq 0}\bigoplus J^tM/J^{t+1}M$$ the associated graded module over the graded ring$$Gr_{J} (R) = \underset{t\geq 0}\bigoplus\ J^t/J^{t+1} .$$ Let  $x_1 , x_2 , \ldots , x_s$ be a system of generators of $J$ and  $(R/J)[T_1 , \ldots ,T_s]$ the polynomial ring of  variables  $T_1 , T_2 , \ldots , T_s$.   The natural epimorphism
$$g: (R/J)[T_1 , \ldots ,T_s] \longrightarrow Gr_{J} (R) $$ leads $K$ to be  an $(R/J)[T_1 , \ldots ,T_s]-$module. We write $K_t= J^tM/J^{t+1}M$ for all $t\geq 0,$ then  $K_0 = M/J M$ is  a noetherian $R/J-$module by the hypothesis. On the other hand, it is easy to check that
$$K_{t+1} = \underset{i=1}{\overset{s}\sum} T_i K_t $$ for all $t\geq 0.$ Thus  $K$
satisfies the conditions  of \cite[1 (i)]{kirart}. Then $K$ is a
noetherian $(R/J )[T_1 , \ldots ,T_s]-$module and so $K$ is a
noetherian $Gr_{J} (R)-$module. Moreover $M$ is $J-$separated by
the hypothesis. Therefore    $M$ is a noetherian $R-$module by
\cite[10.25]{atiint}.  \end{proof}

\begin{theorem}\label{T:dltcnotemddddp} Let $(R,\m)$ be a local noetherian ring and
$M$ a semi-discrete linearly compact  $R-$module.  Then
$H_i^{\frak{m}}(M)$ is a noetherian $\widehat{R}-$module for all
$i\geq 0.$
\end{theorem}
\begin{proof} We prove the theorem by induction on $i$. If  $i=0,$
we have $H^{\frak{m}}_0 (M) \cong \Lambda_{\frak{m}} (M).$ As $M$
is a semi-discrete linearly compact  $R-$module, $M/\m M$ is also
a semi-discrete linearly compact $R/\m-$module. By virtue of
\cite[5.2]{macdua},  $M/\m M$ is a finite dimensional vector
$R/\m-$space. Then $\Lambda_{\m}(M)$ is a noetherian
$\widehat{R}-$module by \cite[7.2.9]{dieele}. Let $i>0.$ Combining
\ref{L:dcddgiaomd} with \ref{L:dktdptdcqdddp} we may assume, by replacing $M$
with $\underset{t>0}\bigcap \m^tM,$ that  there is
an element $x\in \m$ such that $xM=M.$ Then the short exact
sequence of \lc\ modules
$$0\longrightarrow 0:_Mx \longrightarrow M
\overset{x}\longrightarrow M \longrightarrow 0$$ gives rise to a
long exact sequence of local homology modules
$$\ldots\longrightarrow H^{\m}_i(M) \overset{x}\longrightarrow
H^{\m}_i(M) \overset{\delta}\longrightarrow
H^{\m}_{i-1}(0:_Mx)\longrightarrow  \ldots.$$If $0:_Mx =0,$ then
$H^{\m}_i(M)=x H^{\m}_i(M)=\underset{t>0}\bigcap x^tH^{\m}_i(M)=0$
for all $i\geq 0$ by \ref{L:dddpvpksvdn} (ii). We now assume that
$0:_Mx \not=0.$ By the inductive hypothesis, $H^{\m}_{i-1}(0:_Mx)$
is a noetherian $\widehat{R}-$module. Set $H=H^{\m}_i(M),$ we have
$H/xH \cong \Im\delta \subseteq H^{\m}_{i-1}(0:_Mx).$ It follows
that $H/xH$ is a noetherian $\widehat{R}-$module. Thus $H/\hat{\m}
H$ is also a noetherian $\widehat{R}-$module. Moreover,
$\underset{t>0}\bigcap {\hat{\m}}^t H = \underset{t>0}\bigcap
{\m}^t H^{\m}_i(M) = 0.$ Therefore $H$ is a noetherian
$\widehat{R}-$module by \ref{L:tcnote}.
\end{proof}

\begin{theorem}\label{T:dlmddddpnotetcd}  Let $(R,\m)$ be a local noetherian ring
and $M$ a semi-discrete linearly compact $R-$module with $\Ndim
M=d$. Then $H^I_d(M)$ is a noetherian module on $\Lambda_I(R)$.
\end{theorem}
\begin{proof} We argue by induction on $d$.  If   $d=0,$    $M$ is
a  finitely generated $R-$module, and so is  $I-$separated. By
\ref{T:dldtmditach}, $H^I_0(M)\cong \Lambda_I(M)\cong M,$
therefore $H^I_0(M)$ is a noetherian $\Lambda_I(R)-$module. Let
$d>0.$ From  \ref{L:dcddgiaomd}   we have $H^I_d(M)\cong
H^I_d(\underset{t>0}\cap I^tM).$ If $\Ndim (\underset{t>0}\cap
I^tM) < d,$ then $H^I_d(M) = 0$ by \ref{T:dlttdddpcp} and then
there is nothing to prove. If $\Ndim (\underset{t>0}\cap
I^tM) = d,$  by
\ref{L:dktdptdcqdddp} we may assume,  by replacing $M$ with $\underset{t>0}\cap I^tM,$ that there is an element
$x\in I$ such that $xM=M.$ Then, from the short exact sequence of
\lc\ modules
$$0\longrightarrow 0:_Mx \longrightarrow M
\overset{x}\longrightarrow M \longrightarrow 0$$ we get  an exact
sequence of local homology modules$$H^I_d(M)
\overset{x}\longrightarrow H^I_d(M)
\overset{\delta}\longrightarrow H^I_{d-1}(0:_Mx).$$ Note by
\ref{L:cntmoochiamx} that $\Ndim(0:_Mx)\leqslant d-1.$ If
$\Ndim(0:_Mx)<d-1,$ then $H^I_{d-1}(0:_Mx)=0$ by
\ref{T:dlttdddpcp} and therefore
$$H^I_d(M)=xH^I_d(M)=\underset{t>0}\bigcap x^tH^I_d(M)=0$$
by \ref{L:dddpvpksvdn} (ii). Assume  that
$\Ndim(0:_Mx)=d-1$. It follows by the inductive hypothesis that
$H^I_{d-1}(0:_Mx)$ is a noetherian
$\Lambda_I(R)-$module. On the other hand, we have
$H^I_d(M)/xH^I_d(M) \cong \Im \delta\subseteq H^I_{d-1}(0:_Mx).$
Thus $H^I_d(M)/xH^I_d(M)$ is a noetherian $\Lambda_I(R)-$module.
Therefore $H^I_d(M)/JH^I_d(M)$ is a noetherian
$\Lambda_I(R)-$module, where $J=I\Lambda_I(R)$. Moreover, since
$\underset{t>0}\bigcap J^tH^I_d(M) = \underset{t>0}\bigcap
I^tH^I_d(M) = 0$  and $\Lambda_I(R)$ is complete in $J-$adic
topology, $H^I_d(M)$ is a noetherian $\Lambda_I(R)-$module by
\ref{L:tcnote} as required. \end{proof}
\medskip

\section{Macdonald duality}\label{S:dnmdn}
\medskip

Henceforth $(R,\m)$ will be a  local noetherian ring with the
maximal ideal $\m$. Suppose now that the topology  on  $R$ is the
$\m-$adic topology.
\medskip

Let $M$ be an  $R-$module and $E(R/\m)$ the injective envelope of
$R/\m.$ The module $D(M)=\Hom(M,E(R/\m))$ is called  Matlis dual
of $M.$ If $M$ is a Hausdorff linearly topology $R-$module, then
{\it Macdonald dual}  of $M$ is  defined  by $M^*=Hom(M,E(R/\m))$
the set of continuous homomorphisms of $R-$modules (see \cite[\S
9]{macdua}).  In case $(R,\m)$ is local complete, the topology on
$M^*$ is defined as in \cite[8.1]{macdua}. Moreover, if $M$ is
semi-discrete, then the topology of $M^*$ coincides with that
induced on it as a submodule of $E(R/\m)^M,$ where
$E(R/\m)^M=\underset{x\in M}\prod (E(R/\m))_x,$
$(E(R/\m))_x=E(R/\m)$ for all $x\in M$ (see \cite[8.6]{macdua}).

\begin{lemma}\label{L:dnmdnbmlmdnrr} {\rm(see \cite[5.8]{macdua})}   A  \hlt\
$R-$module $M$ is semi-discrete if and only if $D(M)=M^*.$
\end{lemma}

\begin{lemma}\label{L:dktddcltdas}  {\rm(see \cite[5.7]{macdua})}  Let $M$ be a
\hlt\ $R-$module and $u: M\longrightarrow A^*$ a homomorphism.
Then the following statements  are equivalent:

\centerline{ a) $u$ is continuous, \ \ b) $\ker u$ is open, \ \ c)
$\ker u$ is closed.}
\end{lemma}

A \hlt $R-$module is $\m-${\it primary} if each element of $M$ is
annihilated by a power of $\m.$ A \hlt\ $R-$module $M$ is {\it
linearly discrete} if every $\m-$primary quotient of $M$ is
discrete. It should be noted that  if $M$ is linearly discrete,
then $M$ is semi-discrete. The direct limit of a direct system of
linearly discrete $R-$modules is linearly discrete. If $f: M
\longrightarrow N$ is an epimorphism of Hausdorff linearly
topologized $R-$modules in which $M$ is linearly discrete, then
$f$ is continuous (see \cite[6.2, 6.7, 6.8]{macdua}). Then  I. G.
Macdonald \cite{macdua} established the duality between linearly
discrete and linearly compact modules as follows.

\begin{theorem}\label{T:dnmdncplrrnrrlcp} {\rm( \cite[9.3, 9.12, 9.13]{macdua})} Let
$(R, \m)$ be a complete local noetherian ring.

\noindent (i) If $M$ is linearly compact, then $M^*$ is  linearly
discrete (hence semi-discrete).
If $M$ is semi-discrete, then $M^*$ is  linearly compact.

\noindent (ii) If $M$ is linearly compact or linearly discrete,
then we have a topological isomorphism $\omega
:M\overset{\simeq}\longrightarrow M^{**}$.
\end{theorem}

The following duality theorem between local homology and local cohomology modules
is the main result of this section.

\begin{theorem}\label{T:dldndddpvddddp} (i)  Let  $M$
be an $R-$module. Then for all $i\geq 0,$
$$H^I_i (D(M)) \cong D(H^i_I (M)).$$

\noindent (ii) If $M$  is a  linearly compact $R-$module, then for
all $i\geq 0,$
$$H^I_i(M^*) \cong (H_I^i(M))^*.$$ Moreover, if $(R,\m)$ is a complete local noetherian ring, then
$$H_I^i(M^*) \cong (H^I_i(M))^*.$$

\noindent (iii) If $(R,\m)$ is  a complete local noetherian ring
and $M$  a semi-discrete linearly compact $R-$module, then we have
topological isomorphisms of $R-$modules  for all $i\geq 0,$
$$H_I^i(M^*) \cong (H^I_i(M))^*,$$
$$H^I_i(M^*) \cong (H_I^i(M))^*.$$
\end{theorem}

To prove Theorem \ref{T:dldndddpvddddp} some auxiliary lemmas are
necessary. First, we show that the Macdonald dual functor $(-)^*$
is exact on the category of linearly compact $R-$modules and
continuous homomorphisms.

\begin{lemma}\label{L:dnmdnkhop} Let $$ 0 \lr M \overset{f}\lr N \overset{g}\lr P
\lr 0$$ be a short exact sequence of linearly compact $R-$modules,
in which the homomorphisms $f, g$ are continuous. Then the induced
sequence
$$ 0 \lr P^* \overset{g^*}\lr N^* \overset{f^*}\lr M^* \lr 0$$ is exact.\end{lemma}
\begin{proof}  By \cite[5.5]{macdua} $f$ is an open  mapping, so
replace $M$ by $f(M)$ we may assume that $M$ is a close submodule
of $N.$ Therefore, by \cite[5.9]{macdua}, for any continuous
homomorphism $h: M\longrightarrow E(R/\m)$  there is a continuous
homomorphism $\varphi: N\longrightarrow E(R/\m)$ which extends
$h$.  Thus $f^*$ is surjective.   It is easy to see that $g^*$ is
injective and $\Im g^*\subseteq\Ker f^*.$ So it remains to show
that $\Ker f^*\subseteq\Im g^*.$ Let $\psi\in \Ker f^*,$ we have
$\psi(\Ker g)=\psi(f(M))=0.$ Then $\psi$ induces a homomorphism
$\phi: P\longrightarrow E(R/\m)$ such that $\phi \circ g=\psi.$ It
follows $\ker\phi=g(\ker\psi).$ Since  $\psi$ is continuous,
$\ker\psi$ is open by \ref{L:dktddcltdas}. Moreover, $g$ is open,
so  $\ker\phi$ is also open.  Therefore $\phi$ is continuous by
\ref{L:dktddcltdas}. Thus $\psi\in\Im g^*.$ This finishes the
proof.\end{proof}

Note that submodules and homomorphic images of a semi-discrete module are
also semi-discrete. The following consequence shows that the converse is also true in the category of
linearly compact $R-$modules.

\begin{corollary}\label{C:dkdmdlcpnrr} Let $$ 0 \lr M \overset{f}\lr N \overset{g}\lr P \lr 0
$$ be a short exact sequence of linearly compact
$R-$modules with continuous homomorphisms $f, g$. If
$M$ and $P$ are semi-discrete, then $N$ is also
semi-discrete.\end{corollary} \begin{proof}   It follows from
\ref{L:dnmdnbmlmdnrr} and the hypothesis that $P^*=D(P)$ and
$M^*=D(M).$  We now have a commutative diagram
$$\begin{matrix} 0 \lr& P^* &\overset{g^*}\lr& N^*&\overset{f^*}\lr& M^*& \lr 0\\
& \| && \downarrow j && \| \\
0 \lr & D(P)& \overset{D(g)}\lr& D(N)&\overset{D(f)}\lr& D(M)& \lr
0,\end{matrix}$$ in which $j$ is an inclusion and rows are exact
by \ref{L:dnmdnkhop} and \cite[3.16]{rotani}. It follows that
$N^*=D(N),$ thus $N$ is semi-discrete by \ref{L:dnmdnbmlmdnrr}
(i). \end{proof}

\begin{lemma}\label{L:dnmdnctovaext}   Let $N$ be a finitely generated $R-$module
and  $M$  a   linearly compact $R-$module. Then
$$ ({\Tor}_i^R(N,M))^* \cong{\Ext}^i_R(N,M^*), $$
$${\Tor}_i^R(N,M^*) \cong ({\Ext}^i_R(N,M))^*$$for all $i\geq 0.$\end{lemma}
\begin{proof} Let  $${\bold F_\bullet}: \cdots  \longrightarrow F_i
\longrightarrow F_{i-1} \longrightarrow \cdots \longrightarrow F_1
\longrightarrow F_0 \longrightarrow N \longrightarrow 0$$ be a
free resolution of $N,$ in which the free $R-$modules $F_i$ are
finitely generated. Consider  ${\bold F_\bullet}\bigotimes_R M$ as
a complex  of linearly compact $R-$modules with continuous
differentials. Since the Macdonald dual functor $(-)^*$ is exact
on the category of linearly compact $R-$modules and the continuous
homomorphisms by \ref{L:dnmdnkhop},  it follows by \cite[6.1
Theorem 1]{norani} that
$$(H_i({\bold F_\bullet}\otimes_R M))^*\cong H^i(({\bold F_\bullet}\otimes_R M)^*).$$
On the other hand,  by virtue of  \cite[2.5]{macdua} we have
$$({\bold F_\bullet}\otimes_R M)^*\cong {\Hom}_R({\bold
F_\bullet}, M^*).$$ Therefore
\begin{align*} ({\Tor}_i^R(N,M))^*&\cong (H_i({\bold F_\bullet}\otimes_R M))^*\\
&\cong
H^i({\Hom}_R({\bold F_\bullet}, M^*))\\
&\cong {\Ext}^i_R(N,M^*).\end{align*}The proof of the second
isomorphism is similar. \end{proof}

If $M$ is a linearly topologized $R-$module, then the module
$\Ext_R^i(R/I^t,M)$ is also  a linearly topologized $R-$module by
the topology defined as in \ref{L:extcptt}. Since the local
cohomology module $H_I^i(M)=\dlim \Ext^i_R(R/I^t , M)$ is a
quotient module of $\underset{t}{\oplus}  \Ext^i_R(R/I^t , M)$, it
becomes a linearly topologized $R-$module  with the quotient
topology.

\begin{lemma}\label{L:ddddpcmdrrtt}  Let $(R,\m)$ be a complete local noetherian ring. If $M$ is a
semi-discrete linearly compact $R-$module, then $M$ is linearly discrete and therefore the local cohomology
modules $H^i_I(M)$ are linearly discrete $R-$modules for all $i\geq
0.$
\end{lemma}
\begin{proof} We first show that if $M$ is a semi-discrete linearly
compact $R-$module, then $M$ is linearly discrete. Indeed, since
$M$ is semi-discrete, $M^*$ is linearly compact  and  hence
$M^{**}$ is linearly discrete   by \ref{T:dnmdncplrrnrrlcp} (i).
On the other hand, since $M$ is a  linearly compact $R-$module, we
have by \ref{T:dnmdncplrrnrrlcp} (ii) a topological isomorphism
$M\cong M^{**}$. Therefore $M$ is linearly discrete. Now, by the
same argument as in the proof of \ref{L:extcptt} we can prove that
$\{\Ext^i_R(R/I^t , M)\}_t$ is a direct system of semi-discrete
linearly compact $R-$modules with the continuous homomorphisms,
and therefore it is a direct system of linearly discrete modules.
Thus, by \cite[6.7]{macdua} $ H^i_I(M)=\dlim \Ext^i_R(R/I^t , M)$
are linearly discrete for all $i\geq 0$.
\end{proof}

Now we are able to prove the duality theorem
\ref{T:dldndddpvddddp}.

\begin{proof}[Proof of Theorem \ref{T:dldndddpvddddp}] (i) was proved in
\cite[3.3 (ii)]{cuothe}.

\noindent (ii)  Note  by \cite[2.6]{macdua} that for a direct
system $\{M_t\}$ of \hlt\ $R-$modules with the continuous
homomorphisms we have an isomorphism $\inlim M_t^*\cong (\dlim
M_t)^*.$ Moreover, since  $\big\{\Ext_R^i(R/I^t,M)\big\}_t$ forms
a direct system of linearly compact $R-$modules with continuous
homomorphisms by \ref{L:extcptt},  we get by \ref{L:dnmdnctovaext}
that

\begin{align*} H^I_i (M^*) &= \underset{t}{\underleftarrow{\lim}}{\Tor}^R_i (R/I^t , M^*)\\
&\cong \underset{t}{\underleftarrow{\lim}} ({\Ext}_R^i (R/I^t , M))^* \\
&\cong (\underset{t}{\underrightarrow{\lim}} {\Ext}_R^i (R/I^t ,
M))^* = (H^i_I (M))^*. \end{align*} To prove the second
isomorphism note by \cite[9.14]{macdua} that for an inverse system
$\{M_t \}$ of linearly compact modules over complete local
noetherian ring with continuous homomorphisms  we have an
isomorphism $(\inlim M_t)^*\cong \dlim M^*_t$, and that
$\{\Tor^R_i(R/I^t,M)\}_t$ forms an inverse system of linearly
compact $R-$modules with continuous homomorphisms by
\ref{L:torcptt}. It follows by \ref{L:dnmdnctovaext} that
\begin{align*} H^i_I(M^*) &= \dlim {\Ext}^i_R(R/I^t,M^*)\\
&\cong \dlim ({\Tor}^R_i(R/I^t,M))^*\\
&\cong (\inlim {\Tor}^R_i(R/I^t,M))^* = (H^I_i(M))^*. \end{align*}

\noindent (iii) Let us  prove the first isomorphism. From (ii), it
is the  algebraic isomorphism. Thus,  by \cite[6.8]{macdua}, we
only need to show that both $H_I^i(M^*)$ and $(H_i^I(M))^*$ are
linearly discrete. Indeed, it follows from  \ref{P:dddpcpttlcptt}
and \ref{T:dnmdncplrrnrrlcp} (i) that $(H^I_i(M))^*$ is linearly
discrete.  On the other hand, since $M$ is semi-discrete linearly
compact, $M^*$ is linearly compact and linearly discrete.
Therefore the local cohomology modules $H_I^i(M^*)$ are linearly
discrete by \ref{L:ddddpcmdrrtt}, and the first topological
isomorphism is proved. The second topological isomorphism follows
from the first one  and \ref{T:dnmdncplrrnrrlcp} (ii).
\end{proof}

\begin{corollary}\label{C:hqdldddddp} Let $(R,\m)$ be a complete local
noetherian ring.

\noindent (i) If  $M$ is \lc\   $R-$module,  then for  all $i\geq
0,$
$$H_I^i(M) \cong (H^I_i(M^*))^*,$$
$$H^I_i(M) \cong (H_I^i(M^*))^*.$$

\noindent (ii) If $M$ is a semi-discrete linearly compact $R-$module,
then we have  topological isomorphisms of $R-$modules  for all $i\geq 0,$
$$H_I^i(M) \cong (H^I_i(M^*))^*,$$
$$H^I_i(M) \cong (H_I^i(M^*))^*.$$

\end{corollary}
\begin{proof}  (i) follows from \ref{T:dldndddpvddddp} (ii),  \ref{T:dnmdncplrrnrrlcp} (ii).

\noindent (ii) follows from \ref{T:dldndddpvddddp} (iii) and
\ref{T:dnmdncplrrnrrlcp} (ii).
\end{proof}

\bigskip

\section{Local cohomology of semi-discrete\\ linearly compact
modules}\label{S:ddddpmdrr}
\medskip

In this section  $(R,\m)$ is a  local noetherian ring with the
$\m-$adic topology. We denote by $(\widehat R, \hat {\m})$ the
$\m-$adic completion of  $R$ with the  maximal ideal $\hat{ \m}$ and
$\widehat M$ the $\m-$adic completion of the module $M$. Recall
that an artinian $R-$module $A$  has a natural structure as a
module over $\widehat{R}$ as follows (see \cite[1.11]{shaame}):
Let $\hat a= (a_n) \in \widehat R$ and $x\in A$; since $\m^k x=0$
for some positive integer $k$, $a_nx$ is constant for all  large
$n$, and we define $\hat a x$ to be this constant value. Then we
have the following generalization of this fact for \lc\
$R-$modules.

\begin{lemma}\label{L:mdcpnrrtvrm}   Let  $M$ be a \lc\ $R-$module. Then the
following statements are true.

\noindent (i)   $M$ has a natural structure as a \lc\ module over
$\widehat{R}.$ Moreover, a subset $N$ of $M$  is a \lc\
$\widehat{R}-$submodule  if and only if $N$ is a closed
$R-$submodule.

\noindent (ii) Assume in addition that  $M$  is a semi-discrete $R-$module. Then  $M$ is also a semi-discrete  \lc\ $\widehat{R}-$module.
\end{lemma}
\begin{proof} (i)  Assume that $\big\{ U_i\big\}_{i\in J}$ is a
nuclear base of $M$ consisting of submodules. Then  $M \cong
\underset{i\in J}{\underleftarrow{\lim}} M/U_i,$ in which $M/U_i$
is an artinian $R-$module for all $i\in J$  by \cite[3.11, 4.1,
5.5]{macdua}. It  should be noted by \cite[1.11]{shaame} that an
artinian module over a local noetherian ring $(R,\m)$ has a
natural structure as an artinian module over $\widehat{R}$ so that
a subset of $M$ is an $R-$submodule if and only if it is an
$\widehat{R}-$submodules. Thus $\big\{ M/U_i \big\}_{i\in J}$ can
be regard as an inverse system of artinian $\widehat{R}-$modules
with $\widehat{R}-$homomorphisms. Therefore, pass to the inverse
limits, $M$  has a natural structure as a \lc\ module over
$\widehat{R}.$

It is clear that a \lc\ $\widehat R-$submodule of $M$ is a closed
$R-$submodule. Now, if $N$ is a closed $R-$module of $M,$ then $N
\cong\underset{i\in J}{\underleftarrow{\lim}} N/(N\cap U_i).$
Since
$$N/(N\cap U_i)\cong (N+U_i)/U_i\subseteq M/U_i,$$  $N/(N\cap U_i)$
can be considered as an artinian $R-$submodule of $M/U_i,$ so it
is an artinian   $\widehat{R}-$submodule. Moreover, the
homomorphisms of the inverse system  $\big\{ N/(N\cap U_i)
\big\}_{i\in J}$ are induced from the inverse system $\big\{ M/U_i
\big\}_{i\in J}.$ Therefore, by \ref{L:tcmdcptt4} (iv) $N$ is an
\lc\ $\widehat{R}-$submodule of $M$.

\noindent (ii) follows immediately from (i) by the fact that all
submodules of a semi-discrete \lc\ module are closed. \end{proof}

Remember that the
{\it (Krull) dimension}  $\dim_R M$  of a non-zero
$R-$module $M$ is the supremum of lengths of chains of primes in
the support of $M$ if this supremum exists, and $\infty$
otherwise. If $M$ is finitely generated, then $\dim M = \max\{\dim
R/\p\mid \p\in \Ass M\}$. For convenience, we set $\dim M= -1$ if $M=0$.

\begin{corollary}\label{C:cnotevkrtrenrmu}  Let $M$  be a semi-discrete \lc\
$R-$module. Then

\noindent (i) $\Ndim_R M = \Ndim_{\widehat{R}} M;$

\noindent (ii) $\dim_R M = \dim_{\widehat{R}} M.$
\end{corollary}
\begin{proof} (i) follows immediately from \ref{L:mdcpnrrtvrm} (ii) and the
definition of Noetherian dimension.

\noindent (ii) In the special case $M$ is a finitely generated
$R-$module, from \ref{T:dldtmditach} we have $M \cong
\Lambda_{\m}(M) = \widehat{M}$,  and therefore $\dim_R M =
\dim_{\widehat{R}} M$ by [3, 6.1.3]. For any semi-discrete
linearly compact module $M$, there is by \cite[Theorem]{zoslin} a
short exact sequence $0 \longrightarrow N \longrightarrow M
\longrightarrow A \longrightarrow 0$, where $N$ is finitely
generated and $A$ is artinian. As $\dim_R A = \dim_{\widehat{R}} A
=0$ and $\dim_R N = \dim_{\widehat{R}} N,$ we get $\dim_R M =
\dim_{\widehat{R}} M.$
\end{proof}

\begin{remark}\label{R:lytptcpnrrtdtn}  (i) Denote $\mathcal{C}$ the category of semi-discrete
\lc\ $R-$modules. It is well-known that
the category $\mathcal{C}$ contains the category of artinian $R-$modules
and also the category of finitely generated $R-$modules if $R$ is
complete. However, there are many semi-discrete linearly compact
$R-$modules which are neither artinian nor finitely generated. The
first example for this conclusion is the module $K$ in Remark
\ref{R:rmdlttcpvdzg}. More general, let  $R$ be complete ring, $A$
an artinian $R-$module with $\Ndim A>0$  and $N$ a finitely
generated $R-$module with $\dim N>0$. Then $M=A\oplus N$ is
semi-discrete \lc.  Further, let $Q$ be a quotient module of $M,$
then $Q$ is also a semi-discrete \lc\ $R-$module.

\noindent (ii)  If $M\in$ $\mathcal{C}$,  then by
\ref{L:dnmdnbmlmdnrr} the Matlis dual $D(M)$ and the Macdonald dual
$M^*$ are the same. Moreover, $M$  is linearly discrete by
\ref{L:ddddpcmdrrtt}  and can be regarded by  \ref{L:mdcpnrrtvrm}
as an $\widehat R-$module, therefore the Macdonald dual functor
$(-)^*$ is a functor from $\mathcal{C}$ to itself  and we have by
\ref{T:dnmdncplrrnrrlcp} a topological isomorphism $\omega
:M\overset{\simeq}\longrightarrow M^{**}$. Thus, $(-)^{*}$ is an
equivalent functor on the category $\mathcal{C}.$
 \end{remark}

\begin{lemma}\label{L:ndimmsaobangdim}  Let $M$ be a semi-discrete linearly
compact $R-$module. Then $${\Ndim}_R M^* = {\dim}_R M \  \ and \ \
{\Ndim}_R M = {\dim}_R M^*.$$
\end{lemma}
\begin{proof} From \ref{L:mdcpnrrtvrm} (ii) and \ref{C:cnotevkrtrenrmu} we may assume that $(R,\m)$
is a complete ring. If  $M$ is finitely generated $R-$module,
$M^*$ is artinian. Keep in mind  in our case that $M^*=D(M)$, then
the equality $ \Ndim M^* =\dim M$ follows from the well-known
facts of Matlis duality. If $M$ is artinian, then it is clear that
$\Ndim M^* = \dim M=0.$ Suppose now that  $M$ is  semi-discrete
\lc. There is by \cite[Theorem]{zoslin} a short exact sequence
$0\longrightarrow N \longrightarrow M \longrightarrow A
\longrightarrow 0$ in which $N$ is finitely generated and $A$ is
artinian. Thus we get by Macdonald duality an exact sequence
$0\longrightarrow A^* \longrightarrow M^* \longrightarrow N^*
\longrightarrow 0,$ where $N^*$ is artinian and $A^*$ is finitely
generated. Then
\begin{align*}\Ndim M^* &= \max\big\{
\Ndim N^*, \Ndim A^* \big\}\\
& = \max\big\{ \dim N, \dim A \big\} = \dim M .\end{align*} The
second equality follows from \ref{R:lytptcpnrrtdtn},
(ii).\end{proof}

Now we  are able to extend  well-known
results in Grothendieck's local cohomology theory of finitely generated $R-$modules for semi-discrete \lc\ modules.

\begin{theorem}\label{T:mrttddddpg}  Let
$M$ be a non zero semi-discrete linearly compact $R-$module. Then

\noindent (i) $\dim_R M= \max \big\{ i\mid
H_{\m}^i(M)\not=0\big\}$ if   $\dim_R M\not= 1;$

\noindent (ii) $\dim_R (\widehat{M}) = \max \big\{ i\mid
H_{\m}^i(M)\not=0\big\}$ if $\widehat{M}\not= 0.$
\end{theorem}
\begin{proof} (i)  Note by \ref{L:mdcpnrrtvrm} (ii) and \ref{C:cnotevkrtrenrmu} that $M$ is a
semi-discrete \lc\ $\widehat{R}-$module  with $\dim_R M =
\dim_{\widehat{R}} M$. Moreover, the natural homomorphism $f:
R\longrightarrow \widehat{R}$ gives by \cite[4.2.1]{broloc} an
isomorphism $H^d_{\m}(M)\cong H^d_{\widehat{\m}}(M).$ Thus,  we
may assume without any loss of generality that $(R,\m)$ is a
complete local noetherian ring. As $M$ is a semi-discrete linearly
compact $R-$module, $M^*$ is also a semi-discrete linearly compact
$R-$module by \ref{T:dnmdncplrrnrrlcp} (i). Recall by
\cite[5.6]{macdua} that $M^* = 0$ if only if $ M=0$. Then, since
$1\not=  \dim M = \Ndim M^*$, it follows from
\ref{L:ndimmsaobangdim}, \ref{T:dlttdddpcpnrr} (ii) and
\ref{T:dldndddpvddddp} (ii) that
\begin{align*} \dim M = \Ndim M^*
&= \max \big\{ i\mid
H^{\m}_i(M^*)\not=0\big\}\\
&=  \max \big\{ i\mid
H_{\m}^i(M)^*\not=0\big\}\\
&= \max \big\{ i\mid H_{\m}^i(M)\not=0\big\}.\end{align*}

\noindent (ii) The continuous epimorphisms $M\lr M/{\m}^tM$ for
all $t>0$ induce by \ref{L:lnkhmdcptt} a continuous epimorphism
$\pi: M\lr \widehat{M}.$ Moreover $\pi$ is the open homomorphism
by \cite[5.5]{macdua}. Thus $\widehat{M}$ is also a semi-discrete
\lc\ $R-$module. It follows from \ref{C:cnotevkrtrenrmu} that
$\dim_R\widehat{M}=\dim_{\widehat{R}}\widehat{M}.$ Hence, as in
the proof of (i)  we may assume without any loss of generality
that $(R,\m)$ is a complete local noetherian ring. Note that
$H^I_0(M)\cong \Lambda_I(M)$ and $H^0_I(M)\cong \Gamma_I(M),$
hence we have by \ref{T:dldndddpvddddp} (ii)
$$0\not=(\widehat M)^*=(H^{\m}_0(M))^* \cong \Gamma
_{\m}(M^*) .$$ Thus, by virtue of \ref{L:ndimmsaobangdim},
\ref{T:dlttdddpcpnrr} (i) and \ref{T:dldndddpvddddp} (ii) we get
\begin{align*} \dim {\widehat M} = \Ndim (\widehat{M})^* = \Ndim
\Gamma_{\m}( M^*)
 = &\max \big\{ i\mid
H^{\m}_i(M^*) \not=0 \big\}\\
= &\max \big\{ i\mid
(H_{\m}^i(M))^* \not=0 \big\}\\
=&\max \big\{ i\mid H_{\m}^i(M) \not=0 \big\}. \end{align*} The
proof is complete.
\end{proof}

\begin{remark}  \label{R:rmcnt} (i) The condition $\Ndim M\not=1$ in Theorem \ref{T:mrttddddpg}
(i) is necessary. Indeed, take the ring $R$ and the semi-discrete
\lc\ $R-$module $K$ as in Remark \ref{R:rmdlttcpvdzg} and set
$L=K^*$.  It follows from \ref{T:dldndddpvddddp} (iii) and
\ref{R:rmdlttcpvdzg} that
$H_{\m}^i(L) \cong H^{\m}_i(K)^* =0$ for all $i\geq 0$.  Hence
$$\dim L = \Ndim K =1\not= -1 = \max \big\{ i\mid
H_{\m}^i(L)\not=0\big\}.$$

\noindent (ii) The condition $\widehat M\not= 0$ in Theorem
\ref{T:mrttddddpg} (ii)  can also not be dropped as the following
example shows. Set $M=L\oplus A,$ where $A$ is an artinian
$R-$module satisfying $\Width_{\m} A \geq 1$. Then there is an
element $x\in \m$ such that $xM=M$, and therefore $\widehat M=0$.
It is easy to see that $H_{\m}^i(M)=H_{\m}^i(L)=0$ for all $i\geq
1$ and $H_{\m}^0(M)\cong H_{\m}^0(A)=A$. Thus
$$\dim \widehat M =  -1 \not= 0= \max \big\{ i\mid
H_{\m}^i(M)\not=0\big\}.$$
\end{remark}

To complete the unusual behaviour on the vanishing theorem of local cohomology for semi-discrete \lc\ modules we give a characterization of semi-discrete \lc\ modules, whose all local cohomology modules are vanished.

\begin{corollary} \label{C:dtmodddptttb} Let $M$ be a semi-discrete \lc\  $R-$module.
Then $H_{\m}^i(M)=0$ for all $i\geq 0$ if and only if there exists
an element $x\in \m$ such that $xM=M$ and $0:_Mx=0$.
\end{corollary}
\begin{proof} The conclusion follows from \ref{C:hqddpttmi}  by using
\ref{T:dldndddpvddddp} (ii) and \ref{R:lytptcpnrrtdtn} (ii).
\end{proof}

Recall that a sequence of elements $x_1 ,\ldots ,x_r$ in $R$
is said to be an {\it $M-$regular} sequence if $M/(x_1 , \ldots ,
x_r )M \not= 0$ and $M/(x_1 , \ldots , x_{i-1} )M
\overset{x_i}\longrightarrow M/(x_1 , \ldots , x_{i-1} )M$ is
injective for $i=1,\ldots , r.$ Denote by $\depth_I(M)$ the
supremum of the lengths of all  maximal $M-$regular sequences in
$I.$ Then we have

\begin{theorem} \label{T:depthdddpsrcp} Let   $M$ be a \sd\ \lc\
$R-$module such that $M/IM\not=0.$ Then
$${\depth}_I(M)=\inf\{i/H_I^i(M)\not=0\}.$$
\end{theorem}
\begin{proof} Note by \cite[\S
3]{ooimat}  that $x_1 ,\ldots ,x_r$ is an $M-$regular sequence if
and only if it is a $D(M)-$coregular sequence.  Since $M$ is \sd \
\lc, $D(M)=M^*$ and therefore ${\depth}_I(M)={\Width}_I(M^*).$  On
the other hand,  $M^*$ is  a \sd\ \lc\ $R-$module by
\ref{T:dnmdncplrrnrrlcp}(i) and
 $0:_{M^*}I \cong (M/IM)^*\not=0.$
Thus the conclusion follows by virtue of  \ref{T:dtdrdddpcpnrr} and
\ref{T:dldndddpvddddp} (ii).
\end{proof}
Following is the artinianness of local cohomology modules.

\begin{theorem}\label{T:ddddpcpnrrlartin}  Let  $M$ be a  semi-discrete linearly
compact $R-$module with $\dim_R M = d.$ Then  the following
statements are true.

\noindent (i)  The local cohomology modules $H^i_{\m}(M)$ are
artinian $R-$modules for all $i\geq 0;$

\noindent (ii) The  local cohomology module $H_I^d(M)$ is artinian.
\end{theorem}

\begin{proof} Note first that if $A$ is an artinian $\widehat
R-$module, then $A$ is   an artinian $R-$module. Therefore, from
the independent of the base ring of local cohomology and
\ref{L:mdcpnrrtvrm} we may assume without loss of generality that
$R$ is complete. Then, by applying the duality between local
homology and local cohomology \ref{T:dldndddpvddddp}, the
statement (i) follows from \ref{T:dltcnotemddddp} and the
statement (ii) from \ref{T:dlmddddpnotetcd}. \end{proof}

Finally, as an immediate consequence of Theorem
\ref{T:ddddpcpnrrlartin} we get the following well-known result.

\begin{corollary}\label{C:ddddpmdhhsatin}  {\rm (see \cite[7.1.3, 7.1.6]{broloc})}
Let  $M$ be a finitely generated $R-$module with $\dim_R M = d.$
Then the local cohomology modules $H^i_{\m}(M)$ and $H^d_I(M)$ are
artinian $R-$modules for all $i\geq 0.$\end{corollary}

\medskip
\noindent{\bf Acknowledgments.} The  authors have greatly enjoyed
the perusal of the work by I. G. Macdonald \cite{macdua}. Some of
the technical ideals exhibited  in this paper are derived from
this work. They would like to thank professor H. Z\"oschinger for
showing  the module $K$ in Remark \ref{R:rmdlttcpvdzg}. The
authors acknowledge support by the National Basis Research Program
in Natural Science of Vietnam and the Abdus Salam International
Centre for Theoretical Physics, Trieste, Italy.

\end{document}